\documentclass[12pt]{amsart} 
\usepackage{amssymb,amsmath,epsf,colordvi} 
\usepackage{color}
\usepackage{epic}
\usepackage{eepic}
\usepackage{graphicx}

\newtheorem{theorem}[equation]{Theorem}
\newtheorem{lemma}[equation]{Lemma}

\newtheorem{definition}[equation]{Definition}

\title[\includegraphics{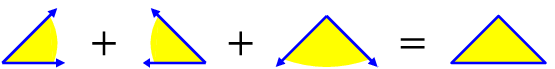}]{
$\raisebox{2pt}{\includegraphics[height=29pt]{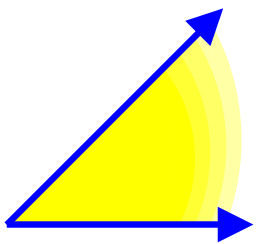}}
 \ \raisebox{9pt}{+}\ \raisebox{2pt}{\includegraphics[height=29pt]{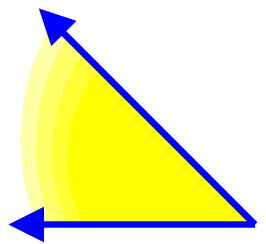}}
 \ \raisebox{9pt}{+}\ \includegraphics[height=27pt]{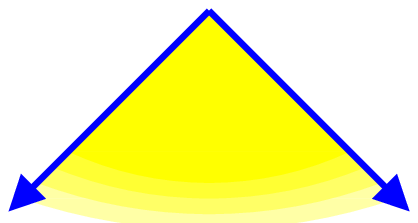} 
 \ \raisebox{9pt}{=}\ \raisebox{3pt}{\includegraphics[height=24pt]{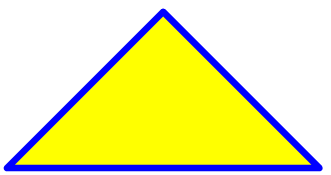}}$\\
 (Formulas of Brion, Lawrence, and Varchenko on rational generating functions
 for cones) 
} 

\author{Matthias Beck}
\address{Department of Mathematics\\
         San Francisco State University\\
         San Francisco, CA 94132\\
         USA}
\email{beck@math.sfsu.org}
\urladdr{http://math.sfsu.edu/beck}

\author{Christian Haase}
\address{Fachbereich Mathematik \& Informatik \\
  Freie Universit\"at Berlin \\
  14195 Berlin\\
  Germany}
\email{christian.haase@math.fu-berlin.de}
\urladdr{http://ehrhart.math.fu-berlin.de}
\author{Frank Sottile}
\address{Department of Mathematics\\
         Texas A\&M University\\
         College Station\\
         TX \ 77843\\
         USA}
\email{sottile@math.tamu.edu}
\urladdr{http://www.math.tamu.edu/\~{}sottile}
\newcommand{\calF}{{\mathcal F}}

\newcommand{\calK}{{\mathcal K}}

\newcommand{\calP}{{\mathcal P}}
\newcommand{\calQ}{{\mathcal Q}}

\newcommand{\V}{\mathcal{V}}

\newcommand{\m}{\mathbf{m}}
\newcommand{\p}{\mathbf{p}}
\newcommand{\s}{\mathbf{s}}
\renewcommand{\v}{\mathbf{v}}
\newcommand{\w}{\mathbf{w}}

\newcommand{\C}{{\mathbb C}}

\newcommand{\N}{{\mathbb N}}
\newcommand{\R}{{\mathbb R}}
\newcommand{\Z}{{\mathbb Z}}

\newcommand{\PL}{\operatorname{PL}}

\newcommand{\conv}{\operatorname{conv}}
\newcommand{\codim}{\operatorname{codim}}
\newcommand{\relint}{\operatorname{relint}}

\newcommand{\QED}{\qquad\epsffile{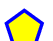}}
\renewcommand{\QED}{\qquad$\Box$}

\newcounter{FNC}[page]
\def\newfootnote#1{{\addtocounter{FNC}{2}$^\fnsymbol{FNC}$%
     \let\thefootnote\relax\footnotetext{$^\fnsymbol{FNC}$#1}}}
\begin{document}

%\begin{abstract}
% We discuss and give elementary proofs of results of Brion and of
% Lawrence--Varchenko on the lattice-point enumerator generating functions for
% polytopes and cones.
% This largely expository note contains a new proof of
% Brion's Formula using irrational decompositions, and a generalization
% of the Lawrence--Varchenko formula.
%\end{abstract}
 
\maketitle
%\tableofcontents

%%%%%%%%%%%%%%%%%%%%%%%%%%%%%%%%%%%%%%%%%%%%%%%%%%%%%%%%%%%%%%%%%%%

%\section*{Introduction}

Our aim is to illustrate two gems of discrete geometry,
namely formulas of Michel Brion~\cite{Br88} and of 
James Lawrence~\cite{La91} and Alexander N.~Varchenko~\cite{varchenko}, which at
first sight seem hard to believe, and which---even after some years of
studying them---still provoke a slight feeling of mystery in us.
Let us start with some examples.

Suppose we would like to list all positive integers. Although there
are many, we may list them compactly in the
form of a generating function:
 \begin{equation}\label{1dray1}
   x^1 + x^2 + x^3 + \cdots\ =\ \sum_{ k>0 } x^k\ =\ \frac{ x }{ 1-x } \ .
 \end{equation}
Let us list, in a similar way, all integers less than or equal to $5$:
 \begin{equation}\label{1dray2}
  \cdots + x^{ -1 } + x^0 + x^1 + x^2 + x^3 + x^4 + x^5\ 
   =\ \sum_{ k \le 5 } x^k\ 
   =\ \frac{ x^5 }{ 1-x^{-1} } \ .
 \end{equation}
Adding the two rational function right-hand sides leads to a miraculous
cancellation 
 \begin{equation}\label{1dpolytope}
  \frac{ x }{ 1-x } + \frac{ x^5 }{ 1-x^{-1} }\ 
  =\ \frac{ x }{ 1-x } + \frac{ x^6 }{ x-1 }\ 
  =\ \frac{ x - x^6 }{ 1-x }\ 
  =\ x + x^2 + x^3 + x^4 + x^5 \,.
 \end{equation}
This sum of rational functions representing two \emph{infinite} series
collapses into a polynomial representing a \emph{finite}
series. 
This is a one-dimensional instance of a theorem
due to Michel Brion. We can think of \eqref{1dray1} as a function
listing the integer points in the ray $[1,\infty)$ and of
\eqref{1dray2} as a function listing the integer points in the ray
$(-\infty,5]$. The respective rational generating functions add up to the
polynomial \eqref{1dpolytope} that lists the integer points in the
interval $[1,5]$. 
Here is a picture of this arithmetic.
 \[
  \begin{picture}(320,65)
   \put(0,15){\includegraphics{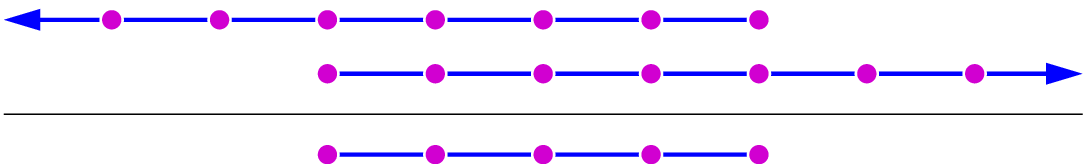}}
%  fig2dev -Leps -m0.4315 Interval.fig Interval.eps
   \put(73,38){$+$}   \put(73,15){$=$}
   \put(92,0){$1$} \put(123,0){$2$}  \put(154,0){$3$}
   \put(185,0){$4$} \put(216,0){$5$} 
  \end{picture}
 \]

Let us move up one dimension. Consider the quadrilateral $\calQ$
with vertices $(0,0)$, $(2,0)$, $(4,2)$, and $(0,2)$. 
 \[
  \begin{picture}(176,88)
   \put(24,10){\includegraphics{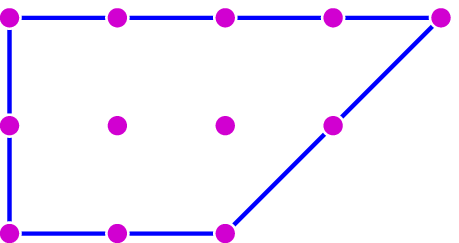}}
%        fig2dev -Leps -m0.4315 Q.fig Q.eps
  % \put(37,57){$\calQ$}  %\put(  8,40){$\calQ$} 
   \put(138,49){$\calQ$} 
   \put(-1,80){$(0,2)$}   \put(155,80){$(4,2)$}
   \put( 0, 0){$(0,0)$}   \put( 90, 0){$(2,0)$}
   \end{picture}
 \]
The analog of
the generating functions \eqref{1dray1} and \eqref{1dray2} are 
the generating functions of the cones at each vertex generated by the
edges at that vertex.
For example, the two edges
touching the origin generate the nonnegative quadrant, which has the generating
function
\[
  \sum_{ m, n \ge 0 } x^m y^n\ =\ 
  \sum_{ m \ge 0 } x^m  \;\cdot\;  \sum_{ n \ge 0 } y^n
  \ =\ \frac{ 1 }{ (1-x) }\cdot \frac{ 1 }{ (1-y) } \ .
\]
The two edges incident to $(0,2)$ generate the cone 
$(0,2) + \R_{ \ge 0 } (0,-2) + \R_{ \ge 0 } (4,0)$, with the generating
function 
\[
  \sum_{ m \ge 0 , n \le 2 } x^m y^n = \frac{ y^2 }{ (1-x) (1-y^{-1}) } \ .
\]
The third such \emph{vertex cone}, at $(4,2)$, is 
$(4,2) + \R_{ \ge 0 } (-4,0) + \R_{ \ge 0 } (-2,-2)$,
which has the generating function
\[
  \frac{ x^4 y^2 }{  (1-x^{-1}) (1-x^{-1}y^{-1})} \ .
\]
Finally, the fourth vertex cone is 
$(2,0) + \R_{ \ge 0 } (2,2) + \R_{ \ge 0 } (-2,0)$, with the generating function
\[
  \frac{ x^2 }{ (1-xy) (1-x^{-1}) } \ .
\]
Inspired by our one-dimensional example above, we add those four
rational functions:
 \begin{align*}
   &\frac{ 1 }{ (1-x)(1-y) } + \frac{ y^2 }{ (1-x)(1-y^{-1}) } + 
    \frac{ x^4 y^2 }{ (1-x^{-1}) (1-x^{-1}y^{-1}) } + \frac{ x^2 }{ (1-xy) (1-x^{-1}) }
   \\
    =& \makebox[1em][l]{}\makebox[1.1em][l]{$y^2$} 
       \makebox[2.5em][l]{$+\, x y^2$}\makebox[2.7em][l]{$+\,x^2 y^2$} + x^2 y^2 + x^4 y^2  \rule{0pt}{14pt}\\
     & \makebox[1em][l]{$+$}\makebox[1.1em][l]{$y$} 
       \makebox[2.5em][l]{$+\, x y$}\makebox[2.7em][l]{$+\,x^2 y $}  + x^3 y \\
     & \makebox[1em][l]{$+$}\makebox[1.1em][l]{$1$} 
       \makebox[2.5em][l]{$+\, x$}\makebox[2.7em][l]{$+\,x^2 $}.
 \end{align*}
The sum of rational functions again collapses to a polynomial, which
encodes precisely those integer points that are contained in
the quadrilateral $\calQ$.

Brion's Theorem says that this magic happens for any polytope $\calP$
in any dimension $d$, provided that $\calP$ has rational vertices.
(More precisely, the edges of $\calP$ have rational directions.) 
The vertex cone $\calK_{\v}$ at vertex $\v$ is the cone with apex $\v$ and
generators the edge directions emanating from $\v$. 
The generating function 
 \[
   \sigma_{ \calK_{\v} } (x)\ :=\ \sum_{ \m \in \calK_{\v} \cap \Z^d } x^\m
 \]
for such a cone is a rational function (again, provided that $\calP$ has
rational vertices).
Here we abbreviate $x^\m$ for $x_1^{ m_1 }
x_2^{ m_2 } \cdots x_d^{ m_d }$. Brion's Formula says that the
rational functions representing the integer points in each vertex
cone sum up to the polynomial $\sigma_\calP(x)$ encoding the integer points in
$\calP$: 
%
% \begin{equation}\label{E:brion}
\[
  \sigma_\calP(x)\ =\ \sum_{\v\textrm{ a vertex of }\calP}
  \sigma_{\calK_{\v}}(x)\ .
 \]% \end{equation}

A second theorem, which shows a similar collapse of 
generating functions of cones, is due (independently) to James Lawrence and
to Alexander Varchenko. 
We illustrate it with the example of the quadrilateral $\calQ$. 
Choose a direction vector $\xi$ that is not perpendicular to any 
edge of $\calQ$, for example we could take $\xi = (2,1)$. 
Now at each vertex $\v$ of $\calQ$, we form a (not necessarily closed) cone
generated by the edge directions $m$ as follows.
If $\w\cdot\xi>0$, then we take its nonnegative span, and if $\w\cdot\xi<0$, we
take its negative span.
\[
  \begin{picture}(370,122)
    \put(0,0){\includegraphics{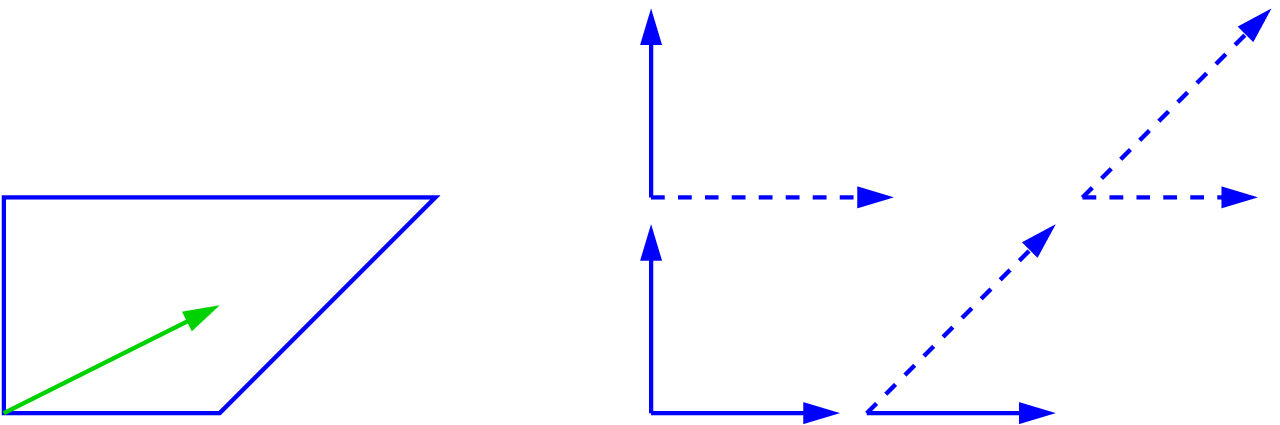}}
%  fig2dev -Leps -m0.4315 Lawrence.fig Lawrence.eps
    \put(40,35){$\xi$}
    \put(110,40){$\calQ$}
   \end{picture}
 \]

For example, the edge directions at the origin are along the positive axes and so this cone
is again the nonnegative quadrant. 
At the vertex $(2,0)$ the edge directions are $(-2,0)$ and $(2,2)$.
The first has negative dot product with $\xi$ and the second has positive dot product,
and so we obtain the half-open cone 
$(2,0) + \R_{< 0} (-2,0) + \R_{\ge 0} (2,2) = (2,0) + \R_{> 0} (2,0) + \R_{\ge 0} (2,2)$. 
At the vertex $(4,2)$ both edge directions have negative dot product with $\xi$ and we
get the open cone $(4,2) + \R_{ > 0 } (0,4) + \R_{ > 0 } (2,2)$, 
and at the vertex $(0,2)$ we get the half-open cone 
$(0,2) + \R_{ \ge 0 } (2,0) + \R_{ > 0 } (0,2)$. 
The respective generating functions are
\[
  \frac{ 1 }{ (1-x)(1-y) } \ , \
  \frac{ x^3 }{  (1-x)(1-xy) } \ , \
  \frac{ x^6 y^3 }{ (1-xy) (1-y) } \ , \  \text{ and } \ 
  \frac{ y^3 }{ (1-x) (1-y) } \ .
\]
Now we add them with signs according to the parity of the number of negative
$(\w\cdot\xi<0)$ edge directions $\w$ at the vertex.  
In our example, we obtain
 \begin{align*}
  &\frac{ 1 }{ (1-x)(1-y) }
   - \frac{ x^3 }{ (1-x) (1-xy) }
   + \frac{ x^6 y^3 }{ (1-xy) (1-y) }
   - \frac{ y^3 }{ (1-x) (1-y) } \\
     =& \makebox[1em][l]{}\makebox[1.1em][l]{$y^2$} 
       \makebox[2.5em][l]{$+\, x y^2$}\makebox[2.7em][l]{$+\,x^2 y^2$} 
            + x^2 y^2 + x^4 y^2  \rule{0pt}{14pt}\\
     & \makebox[1em][l]{$+$}\makebox[1.1em][l]{$y$} 
       \makebox[2.5em][l]{$+\, x y$}\makebox[2.7em][l]{$+\,x^2 y $}  + x^3 y \\
     & \makebox[1em][l]{$+$}\makebox[1.1em][l]{$1$} 
       \makebox[2.5em][l]{$+\, x$}\makebox[2.7em][l]{$+\,x^2 $}.
 \end{align*}
This sum of rational functions again collapses to the polynomial that encodes
the integer points in $\calQ$. 
This should be clear here, for the integer points in the nonnegative quadrant
are counted with a sign $\pm$, depending upon the cone in which they lie, and
these coefficients cancel except for the integer points in the polytope
$\calQ$.

The identity illustrated by this example works for any \emph{simple}
polyope---a $d$-polytope where every vertex meets exactly $d$
edges. 
Given a simple polytope, choose a direction vector $\xi\in\R^d$
that is not perpendicular to any edge direction.
Let $E^+_{\v}(\xi)$ be the edge directions $\w$ at a vertex $\v$ with
$\w\cdot \xi>0$
and $E^-_{\v}(\xi)$ be those with $\w\cdot \xi<0$.
Define the cone
 \[
    \calK_{\xi,\v} :=
    \v\ +\ \sum_{\w\in E^+_{\v}(\xi)} \R_{\geq0} \w\ \ 
      +\ \sum_{\w\in E^-_{\v}(\xi)} \R_{<0} \w\,.
 \]
This is the analogue of the cones in our previous example.
The Lawrence--Varchenko Formula says that adding the rational functions of
these cones with appropriate signs gives the  polynomial
$\sigma_\calP(x)$ encoding the integer points in $\calP$: 
%
% \begin{equation}\label{E:lawrence}
\[
    \sigma_{\calP}(x)\ =\ \sum_{\v\textrm{ a vertex of }\calP}
             (-1)^{|E^-_{\v}(\xi)|} \, \sigma_{\calK_{\xi,\v}}(x)\ .
\]% \end{equation}
Here, $\sigma_{\calK_{\xi,\v}}(x)$ is the generating function encoding the
integer points in the cone $\calK_{\xi,\v}$.
An interesting feature of this identity, which also distinguishes it
from Brion's Formula, is that the power series generating functions have
a common region of convergence.
Also, it holds without any restriction that the polytope be rational.
In the general case, the generating functions of the cones are holomorphic 
functions, which we can add, as they have a common domain (the common
region of convergence).

%
%
%
%%%%%%%%%%%%%%%%%%%%%%%%%%%%%%%%%%%%%%%%%%%%%%%%%%%%%%%%%%%%%%%%%%%%%
\section*{Proofs}
Brion's original proof of his formula~\cite{Br88} used the
Lefschetz--Riemann--Roch theorem in equivariant $K$-theory~\cite{BFQ}
applied to a singular toric variety. 
Fortunately for us, the remarkable formulas of Brion and of Lawrence--Varchenko 
now have easy proofs, based on counting.

\begin{figure}[htb]
\[
%  fig2dev -Leps -m0.2315 P.fig P.eps
%  fig2dev -Leps -m0.2315 PT.fig PT.eps
 \begin{picture}(145,65)(-15,0)
   \put(0,0){\includegraphics{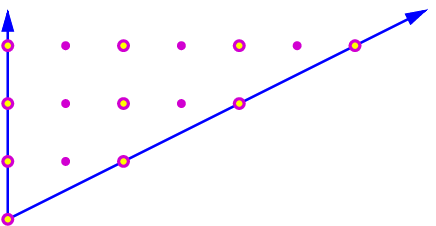}}
   \put(-15,24){$\calK$}
 \end{picture}
  \qquad \qquad
  \begin{picture}(145,65)(-15,0)
   \put(0,0){\includegraphics{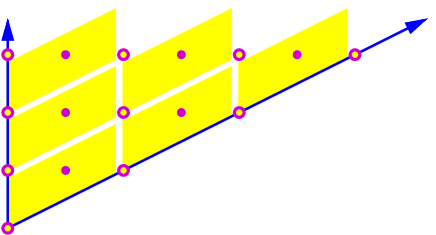}}
   \put(-15,24){$\calK$}
   \put(40,3){$\calP$} \put(37,7){\vector(-4,1){25}}
 \end{picture}
\]
\caption{Tiling a simple cone by translates of its fundamental parallelepiped.}
 \label{F:Tiling}
\end{figure}

Let us first consider an example based on the cone
$\calK=\R_{\geq0}(0,1)+\R_{\geq0}(2,1)$.
The open circles in the picture on the left in Figure~\ref{F:Tiling}
%%%%%%%%%%%%%%%%%%%%%%%%%%%%%%%%%%%%%%%%%%%%%%%%
%%%%%%%%%%%%%%%%%%%%%%%%%%%%%%%%%%%%%%%%%%%%%%%%%%%%%%%%%%%%%%%%%%
represent the semigroup $\N(0,1)+\N(2,1)$, which is a
proper subsemigroup of the integer points $\calK\cap\Z^2$ in $\calK$.
The picture on the right shows how translates of the fundamental half-open
parallelepipied $\calP$ by this subsemigroup cover $\calK$.
This gives the formula
\[
  \sigma_\calK(x)\ =\ \sigma_\calP(x)\cdot
    \sum_{ m,n \ge 0} x^m (x^2y)^n\ =\ 
   \frac{1+xy}{(1-x)(1-x^2y)}\ ,
\]
as the fundamental parallelepiped $\calP$ contains two integer points, the origin
and the point $(1,1)$.

A simple rational cone in $\R^d$ has the form
\[
  \calK\ :=\ \left\{ \v +\sum_{i=1}^d \lambda_i\w_i\mid \lambda_i\in\R_{\geq0}\right\}
       \ =\ \v + \sum_{ i=1 }^d \R_{\geq 0} \w_i\ ,
\]
where $\w_1,\dotsc,\w_d\in\Z^d$ are linearly independent.
This cone is tiled by the 
$(\N\w_1+\dotsb+\N\w_d)$-translates of the half-open parallelepiped
\[
  \calP\ :=\ \left\{\v + \sum_{i=1}^d \lambda_i\w_i\mid 0\leq \lambda_i<1\right\}\ .
\]
The generating function for $\calP$ is the polynomial
\[
   \sigma_{\calP}(x)\ =\ \sum_{\m\in\calP\cap\Z^d} x^{\m}\ ,
\]
and so the generating function for $\calK$ is
 \[
  \sigma_{\calK}(x)\ =\ \sum_{\alpha\in\N\w_1+\dotsb+\N\w_d} 
   x^{\alpha}\cdot \sigma_{\calP}(x)
   \ =\  \frac{\sigma_{\calP}(x)}{(1-x^{\w_1})\dotsb(1-x^{\w_d})}\ ,
 \]
which is a rational function.
This formula and its proof do not require that the apex $\v$ be
rational, but only that the generators $\w_i$ of the cone be linearly
independent vectors in $\Z^d$.

A rational cone $\calK$ with apex $\v$ and generators
$\w_1,\dotsc,\w_n\in\Z^d$  has the form
\[
  \calK\ =\ \v + \R_{\geq0}\w_1+\dotsb+\R_{\geq0}\w_n\ .
\]
If there is a vector $\xi\in\R^d$ with
$\xi\cdot\w_i>0$ for $i=1,\dotsc,n$, then $\calK$ is {\it strictly convex}.
A fundamental result on convexity~\cite[Lemma VIII.2.3]{barvinokBook} 
is that $\calK$ may be decomposed into simple 
cones $\calK_1,\dotsc,\calK_l$ having pairwise disjoint interiors, each with
apex $\v$ and generated by $d$ of the generators $\w_1,\dotsc,\w_n$ of $\calK$.
We would like to add the generating functions for each cone $\calK_i$ to
obtain the generating function for $\calK$.
However, some of the cones may have lattice points in common, and some device is
needed to treat the subsequent overcounting.

An elegant way to do this is to avoid the overcounting altogether
by translating all the cones~\cite{BS05}.
We explain this.
There exists a short vector $\s\in\R^d$ such that 
 \begin{equation}\label{Eq:integer_shift}
   \calK\cap\Z^d\ =\ (\s+\calK)\cap\Z^d\ ,
 \end{equation}
and no facet of any cone $\s+\calK_1,\dotsc,\s+\calK_l$ contains any integer
points. 
This gives the disjoint {\it irrational decomposition}
\[
   \calK\cap\Z^d\ =\ 
   (\s+\calK_1)\cap\Z^d \sqcup\dotsb\sqcup   (\s+\calK_l)\cap\Z^d \ ,
\]
and so
 \begin{equation}\label{Eq:sigma_calK}
   \sigma_{\calK}(x)\ =\ \sum_{\m\in\calK\cap\Z^d} x^{\m}\ =\ 
   \sum_{i=1}^l \sigma_{\s+\calK_i}(x)
 \end{equation}
is a rational function.

For example, suppose that $\calK$ is the cone in $\R^3$ with apex the origin
and generators
\[
   \w_1\ =\ (1,0,1),\quad
   \w_2\ =\ (0,1,1),\quad
   \w_3\ =\ (0,-1,1),\quad\mbox{and}\quad
   \w_4\ =\ (-1,0,1)\,.
\]
If we let $\calK_1$ be the simple cone with generators $\w_1,\w_2,\w_3$ and
$\calK_2$ be the simple cone with 
generators $\w_2,\w_3,\w_4$, then $\calK_1$ and $\calK_2$ decompose $\calK$ into
simple cones.
If $\s=(\frac{1}{8},0,-\frac{1}{3})$, then~\eqref{Eq:integer_shift} holds, and
no facet of $\s+\calK_1$ or of $\s+\calK_2$ contains any integer points.
We display these cones, together with their integer points having
$z$-coordinate 0, 1, or 2.
\[
  \begin{picture}(240,130)(-28,0)
   \put(15,0){\includegraphics[height=130pt]{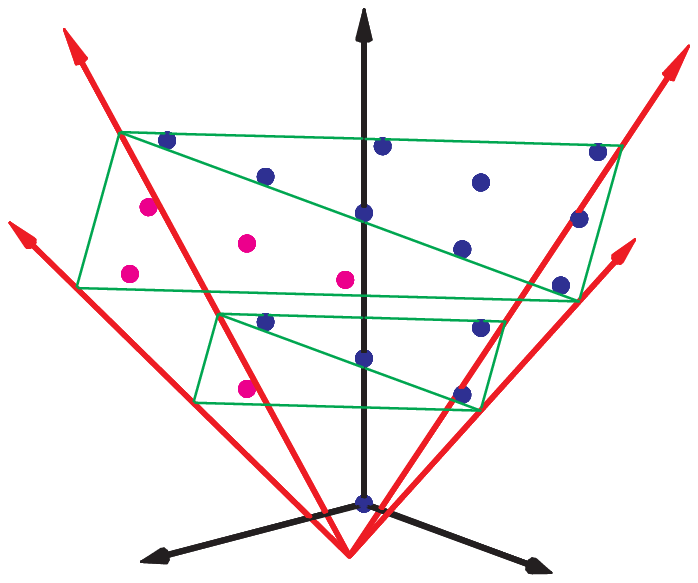}}
   \put(-28,92){$\s+\calK_1$}   \put(9,88){\vector(4,-1){35}}
   \put(177,86){$\s+\calK_2$}\put(173,89){\vector(-1,0){35}}
   \put(135,21){$\s$}\put(133,23){\vector(-2,-1){35}}
   \put(35,117){$\w_3$}  \put(150,117){$\w_4$}
   \put( 5, 68){$\w_1$}  \put(160, 65){$\w_2$}
   \put(41,7){$x$} \put(141,5){$y$} \put(80,120){$z$}
   \end{picture}
\]
The cone $\s+\calK_1$ contains the 5 magenta points shown with positive first
coordinate, while $\s+\calK_2$ contains the other displayed points.
Their integer generating functions are
\begin{eqnarray*}
 \sigma_{\s+\calK_1}(x)&=&\frac{x+xz}{(1-yz)(1-y^{-1}z)(1-xz)}\,,\\
 \sigma_{\s+\calK_2}(x)&=&\frac{1+z}{(1-yz)(1-y^{-1}z)(1-x^{-1}z)}\,,\quad\mbox{and}\\
 \sigma_{\calK}(x)&=&\frac{(1+x)(1-z^2)}{(1-yz)(1-y^{-1}z)(1-xz)(1-x^{-1}z)}\ .
\end{eqnarray*}
Then $\sigma_{\s+\calK_1}(x)+\sigma_{\s+\calK_2}(x)=\sigma_{\calK}(x)$, as
\[
   (x+xz)(1-x^{-1}z) + (1+z)(1-xz)\ =\ 
   1+x-z^2-xz^2\ =\ (1+x)(1-z^2)\,.
\]

While the cones that appear in the Lawrence--Varchenko
formula are all simple, and those in Brion's formula are strictly convex,
we use yet more general cones in their proof.
A rational (closed) halfspace is the convex subset of $\R^d$ defined by 
\[
   \{x\in\R^d\mid \w\cdot x\geq b\}\,,
\]
where $\w\in\Z^d$ and $b\in\R$.
Its boundary is the rational hyperplane $\{x\in\R^d\mid\w\cdot x=b\}$.
A (closed) cone $\calK$ is the interection of finitely many closed halfspaces
whose boundary hyperplanes have some point in common.
We assume this intersection is irredundant.
The {\it apex} of $\calK$ is the intersection of these boundary hyperplanes,
which is an affine subspace.

The generating function for the integer points in $\calK$ is the formal Laurent
series
 \begin{equation}\label{Eq:FGS}
   S_{\calK}\ :=\ \sum_{\m\in\calK} x^{\m}\ .
 \end{equation}
This formal series makes sense as a rational function only if $\calK$ is
strictly convex, that is, if its apex is a single point.
Otherwise,  the apex is a rational affine subspace $L$, and the cone
$\calK$ is stable under translation by any integer vector $\w$ that is parallel
to $L$.
If $\m\in\calK\cap\Z^d$, then the series $S_{\calK}$ contains the series 
\[
   x^{\m}\cdot \sum_{n\in\Z} x^{n\w}
\]
as a subsum.
As this converges only for $x=0$, the series $S_{\calK}$
converges only for $x=0$.  

We relate these formal Laurent series to rational functions.
The product of a formal series and a polynomial is another formal 
series. 
Thus the additive group $\C[[x_1^{\pm1},\dotsc,x^{\pm1}_d]]$ of
formal Laurent series is a module over the ring 
$\C[x_1^{\pm1},\dotsc,x_d^{\pm1}]$ of Laurent polynomials. 
The space $\PL$ of {\em polyhedral Laurent series} is the
$\C[x_1^{\pm1},\dotsc,x_d^{\pm1}]$-submodule of
$\C[[x_1^{\pm1},\dotsc,x^{\pm1}_d]]$ generated by the set of 
formal series  
 \[
    \{S_\calK\mid \calK\textrm{ is a simple rational cone}\}\,.
 \]
Since any rational cone may be triangulated by simple cones, 
$\PL$ contains the integer generating series of all rational cones.  

Let $\C(x_1,\dotsc,x_d)$ be the field of rational functions on $\C^d$, which is
the quotient field of $\C[x_1^{\pm1},\dotsc,x_d^{\pm1}]$. 
According to Ishida~\cite{Is90},
the proof of the following theorem is due to Brion.

\begin{theorem}\label{T:unique}
  There is a unique homomorphism of\/ $\C[x_1^{\pm1},\dotsc,x_d^{\pm1}]$-modules
\[
  \varphi\ \colon\ \PL\ \longrightarrow \C(x_1,\dotsc,x_d)\,,
\]
 such that $\varphi(S_\calK)=\sigma_\calK$ for every simple cone $\calK$ in
 $\R^d$. 
\end{theorem}

\noindent{\it Proof. }
Given a simple rational cone $\calK = v + \langle \w_1,\dots,\w_d\rangle$ with
fundamental parallelepiped $\calP$, we have
\[
  \prod_{i=1}^d (1-x^{\w_i}) \cdot S_\calK\ =\ \sigma_{\calP}(x)\,.
\]
Hence, for each $S\in \PL$, there is a nonzero Laurent
polynomial $g\in\C[x_1^{\pm1},\dotsc,x_d^{\pm1}]$ such that 
$g S = f \in \C[x_1^{\pm1},\dotsc,x_d^{\pm1}]$.
If we define $\varphi(S):= f/g\in \C(x_1,\dotsc,x_d)$, then
$\varphi(S)$ is independent of the choice of $g$.
This defines the required homomorphism.\quad\QED \medskip 

The map $\varphi$ takes care of the nonconvergence of the generating series
$S_{\calK}$ when $\calK$ is not strictly convex.

\begin{lemma}\label{L:Zero}
 If a rational polyhedral cone $\calK$ is not strictly convex, then
 $\varphi(S_\calK)=0$. 
\end{lemma}

\noindent{\it Proof. }
 Let $\calK$ be a rational polyhedral cone that is not strictly
 convex. 
 Then there is a nonzero vector $\w\in\Z^d$ such that 
 $\w+\calK=\calK$, and so $x^{\w}\cdot S_\calK=S_\calK$. 
 Thus $x^{\w} \varphi(S_\calK)=\varphi(S_\calK)$. 
 Since $1-x^{\w}$ is not a zero-divisor in $\C(x_1,\dotsc,x_d)$, we
 conclude that  $\varphi(S_\calK)=0$. 
\QED\medskip

We now establish Brion's Formula, first for a simplex, and then use irrational
decomposition for the general case.
(A $d$-dimensional simplex is the intersection of $d{+}1$ halfspaces, one for each facet.)

For a face $F$ of the simplex $\calP$, let $\calK_F$ be the tangent cone to $F$,
which is the intersection of the halfspaces corresponding to the $d-\dim(F)$
facets containing $F$.
Let $\emptyset$ be the empty face of $\calP$, which has dimension $-1$.
Its tangent cone is $\calP$.

\begin{theorem}\label{T:simplex}
 If $P$ is a simplex, then 
\begin{equation}\label{E:simplex}
    0\ =\ \sum_{F} (-1)^{\dim(F)} S_{\calK_F}\,,
\end{equation}
 the sum over all faces of $P$.
\end{theorem}

\noindent{\it Proof.}
 Consider the coefficient of $x^{\m}$ for some $\m\in\Z^d$ in the sum on the
 right. 
 Then $\m$ lies in the tangent cone $\calK_F$ to a unique face $F$ of
 minimal dimension, as $P$ is a simplex.
 The coefficient of $x^{\m}$ in the sum becomes
\[
  \sum_{G\supseteq F} (-1)^{\dim(G)}\,.
\]
 But this vanishes, as every interval in the face poset of $P$ is a
 Boolean lattice.\quad\QED \medskip 

Now we apply the evaluation map $\varphi$ of Theorem~\ref{T:unique} to the
formula~\eqref{E:simplex}. 
Lemma \ref{L:Zero} implies that $\varphi(S_{\calK_F})= 0$ except when
$F=\emptyset$ or $F$ is a vertex, and then 
$\varphi(S_{\calK_F})=\sigma_{\calK_F}(x)$. 
This gives
\[
    0\ =\ - \sigma_{\calP} (x) + 
    \sum_{ \v \text{ a vertex of } \calP } \sigma_{ \calK_{\v} } (x) \,,
\]
which is Brion's Formula for simplices.

Just as for rational cones, every polytope $\calP$ may be decomposed into simplices
$\calP_1, \dots, \calP_l$ having pairwise disjoint interiors, using only the vertices of $\calP$.
\[
   \calP\ =\ \calP_1 \cup \dotsb\cup \calP_l\ .
\]
Then there exists a small real number $\epsilon>0$ and a short vector $\s$ such that
if we set
\[
  \calP'\ :=\ \s+(1+\epsilon)\calP\quad\mbox{and}\quad
  \calP'_i\ :=\ \s+(1+\epsilon)\calP_i\quad\mbox{for}\ i=1,\dotsc,l\,,
\]
then $\calP'\cap\Z^d=\calP\cap\Z^d$, and no hyperplane supporting any facet of
any simplex $\calP'_i$ meets $\Z^d$.
If we write $\calK(\calQ)_{\w}$ for the tangent cone to a polytope $\calQ$ at a
vertex $\w$, then for $\v$ a vertex of $\calP$ with $\v'=(1+\epsilon)\v+\s$ the
coresponding vertex of $\calP'$, we have 
$\calK(\calP')_{\v'}\cap\Z^d = \calK(\calP)_{\v}\cap\Z^d$ and so this is an irrational
decomposition. 
Then
\begin{eqnarray*}
 \sum_{\v \text{ a vertex of }\calP}\sigma_{\calK(\calP)_{\v}}(x)&=&
 \sum_{\v \text{ a vertex of }\calP'}\sigma_{\calK(\calP')_{\v}}(x)\\
 &=&\sum_{i=1}^l\sum_{\v\text{ a vertex of }\calP'_i }\sigma_{\calK(\calP'_i)_{\v}}(x)\\
 &=&\sum_{i=1}^l \sigma_{\calP_i}(x)\ =\ \sigma_{\calP'}(x)\ =\ \sigma_{\calP}(x)\ .
\end{eqnarray*}
The second equality holds because the vertex cones $\calK(\calP'_i)_{\v}$ form an
irrational decomposition of the vertex cone $\calK(\calP')_{\v}$, and 
because the same is true for the polytopes.
This completes our proof of Brion's Formula.

Consider the quadrilateral $\calQ$, which may be triangulated by adding
an edge between the vertices $(2,0)$ and $(0,2)$.
Let $\epsilon=\frac{1}{4}$ and $\s=(-\frac{1}{2},-\frac{1}{4})$.
Then $(1+\epsilon)\calQ+\s$ has vertices
\[
  ({\textstyle -\frac{1}{2}},{\textstyle -\frac{1}{4}}), \quad
  (2,{\textstyle -\frac{1}{4}}), \quad
  ({\textstyle -\frac{1}{2}},2+{\textstyle \frac{1}{4}}), \quad
  (4+{\textstyle \frac{1}{2}},2+{\textstyle \frac{1}{4}})\ .
\]
We display the resulting irrational decomposition.
\[
  \begin{picture}(130,73)(26,7)
   \put(24,10){\includegraphics{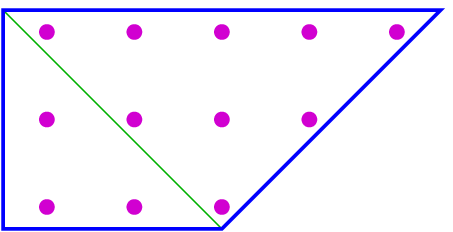}}
%        fig2dev -Leps -m0.35 QT.fig QT.eps
  % \put(37,57){$\calQ$}  %\put(  8,40){$\calQ$} 
   \put(138,46){$\calQ'$} 
   \end{picture}
\]

We use the map $\varphi$ to deduce a very general form of the
Lawrence--Varchenko formula.
Let $\calP$ be a simple polytope, and for each vertex $\v$ of $\calP$ choose a
vector $\xi_{\v}$ that is not perpendicular to any edge direction at $\v$.
Form the cone $\calK_{\xi_{\v}, \v}$ as before.
Then we have
\begin{equation}\label{E:LV}
    \sigma_{\calP}(x)\ =\ \sum_{\v\textrm{ a vertex of }\calP}
             (-1)^{|E^-_{\v}(\xi_{\v})|} \, \sigma_{\calK_{\xi_{\v},\v}}(x)\ .
\end{equation}
Brion's formula is the special case when each vector $\xi_\v$
points into the interior of the polytope.
We establish~\eqref{E:LV} by showing that the sum on the right does not change when
any of the vectors $\xi_{\v}$ are rotated.

Pick a vertex $\v$ and vectors $\xi,\xi'$ that are not perpendicular to any
edge direction at $\v$ such that $\xi\cdot \w$ and $\xi\cdot \w'$ have the same
sign for all except one edge direction $\m$ at $\v$.
Then $\calK_{\xi,\v}$ and $\calK_{\xi',\v}$ are disjoint and their union is the
(possibly) half-open  
cone $\calK$ generated by the edge directions $\w$ at $\v$ such that 
$\xi\cdot\w$ and $\xi'\cdot\w$ have the same sign, but with apex the affine line
$\v+\R\m$.
Thus we have the identity of rational formal series
\[
   S_{\calK_{\xi,\v}} - S_{\calK}\ =\ -S_{\calK_{\xi',\v}}\ .
\]
Applying the evaluation map $\varphi$ gives
\[
  \sigma_{\calK_{\xi,\v}}(x)\ =\ -\sigma_{\calK_{\xi',\v}}(x)\ ,
\]
which proves the claim, and the generalized Lawrence--Varchenko formula~\eqref{E:LV}.

%%%%%%%%%%%%%%%%%%%%%%%%%%%%%%%%%%%%%%%%%%%%%%%%%%%%%%%%%%%%%%%%%%%%%
\section*{Valuations}
Valuations provide a conceptual approach to these
ideas. Once the theory is set up, both Brion's Formula and the 
Lawrence--Varchenko Formula are easy corollaries of
duality being a valuation. 
We are indebted to Sasha Barvinok who pointed out this correspondence
to the second author during a coffee break at the 2005 Park City
Mathematical Institute.
Let us explain.

Consider the vector space of all functions $\R^d \to \R$.
Let $\V$ be the subspace that is generated by indicator
functions of polyhedra:
 \begin{equation*}
  [\calP]\ \colon\ x\ \mapsto\ 
  \begin{cases}
    1 & \text{ if } x \in \calP , \\
    0 & \text{ if } x \not\in \calP .
  \end{cases}
 \end{equation*}
We add these functions point-wise.
For example, if $d=1$, and $\calP=[0,2]$, $\calQ=[1,3]$, then
$[\calP]+[\calQ]$ takes the value $1$ along $[0,1)$ and $(2,3]$, the
value $2$ along $[1,2]$, and vanishes everywhere else.
\[
  \begin{picture}(120,65)(2.5,0)
   \put(2.5,11){\includegraphics{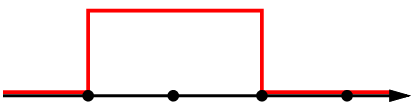}}
%fig2dev -Leps -m0.227 V02.fig V02.eps
   \put(25,0){$0$}    \put( 50,0){$1$} 
   \put(75,0){$2$}    \put(100,0){$3$} 
  \end{picture}
    \quad \raisebox{10mm}{+} \quad 
  \begin{picture}(120,65)(2.5,0)
   \put(2.5,11){\includegraphics{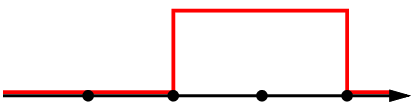}}
%fig2dev -Leps -m0.227 V13.fig V13.eps
   \put(25,0){$0$}    \put( 50,0){$1$} 
   \put(75,0){$2$}    \put(100,0){$3$} 
 \end{picture}
    \quad \raisebox{10mm}{=} \quad 
 \begin{picture}(120,65)(2.5,0)
   \put(2.5,11){\includegraphics{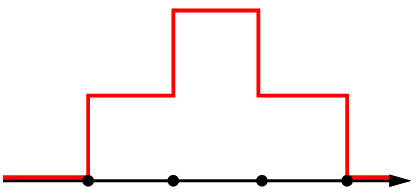}}
%fig2dev -Leps -m0.227 V0123.fig V0123.eps
   \put(25,0){$0$}    \put( 50,0){$1$} 
   \put(75,0){$2$}    \put(100,0){$3$} 
 \end{picture}
\]
Already this simple example shows that our generators do not form a
basis: they are linearly dependent. For $\calP'=[0,3]$ and
$\calQ'=[1,2]$, we get the same sum.
\[
  \begin{picture}(120,65)(2.5,0)
   \put(2.5,11){\includegraphics{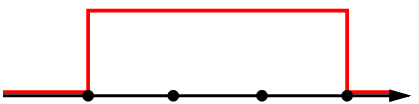}}
%fig2dev -Leps -m0.227 V03.fig V03.eps
   \put(25,0){$0$}    \put( 50,0){$1$} 
   \put(75,0){$2$}    \put(100,0){$3$} 
  \end{picture}
    \quad \raisebox{10mm}{+} \quad 
  \begin{picture}(120,65)(2.5,0)
   \put(2.5,11){\includegraphics{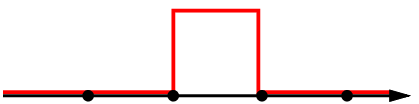}}
%fig2dev -Leps -m0.227 V12.fig V12.eps
   \put(25,0){$0$}    \put( 50,0){$1$} 
   \put(75,0){$2$}    \put(100,0){$3$} 
 \end{picture}
    \quad \raisebox{10mm}{=} \quad 
 \begin{picture}(120,65)(2.5,0)
   \put(2.5,11){\includegraphics{figures/V0123}}
%fig2dev -Leps -m0.227 V0123.fig V0123.eps
   \put(25,0){$0$}    \put( 50,0){$1$} 
   \put(75,0){$2$}    \put(100,0){$3$} 
 \end{picture}
\]
But this is the only thing that can happen.

\begin{theorem}[\cite{Gr78,Vo57}]
  The linear space of relations among the indicator functions $[\calP]$ of
  convex polyhedra is generated by the relations $[\calP]+[\calQ]=[\calP
  \cup \calQ]+[\calP \cap \calQ]$ where $\calP$ and $\calQ$ run over
  polyhedra for which $\calP \cup \calQ$ is convex.
\end{theorem}

A {\it valuation} is a linear map $\nu\colon\V \to V$, where $V$ is some
vector space. 
Some standard examples are
%%%%%%%%%%%%%%%%%%%%%%%%%%%%%%%%%%%%%%%%%%%%%%%%%%%%
 \begin{center}
  \begin{tabular}{|l||c|}\hline
     $V$ & $\nu(\calP)$\rule{0pt}{12pt}\\\hline\hline
   $\R^d$& $\operatorname{vol}(\calP)$\rule{0pt}{12pt}\\\hline
   $\PL$ & $S_\calP(x)$\rule{0pt}{12pt}\\\hline
   $\C(x_1,\dotsc,x_d)$ 
         &$\sigma_\calP(x)$\rule{0pt}{12pt}\\\hline
   $\R^d$& 1\rule{0pt}{12pt}\\\hline
  \end{tabular} \ .
 \end{center}
%%%%%%%%%%%%%%%%%%%%%%%%%%%%%%%%%%%%%%%%%%%%%%%%%%%
That $\sigma_\calP(x)$ is a valuation is a deep result of Khovanskii-Pukhlikov~\cite{KhPu}
and of Lawrence~\cite{La88}.
The last example is called the Euler characteristic.
This valuation is surprisingly useful. For example, it can
be used to prove Theorem~\ref{thm:polarity} below.

The most interesting valuation for us comes from the polar construction.
The {\it polar} $\calP^\vee$ of a polyhedron $\calP$ is the polyhedron given by
 \[
    \calP^\vee\ :=\ \{ x \mid \langle x, y \rangle \le 1 \text{ for all }
    y \in \calP \}\,.
 \]

It is instructive to work through some examples.
\begin{enumerate}
\item
   \begin{tabular}{ccc}
      \includegraphics{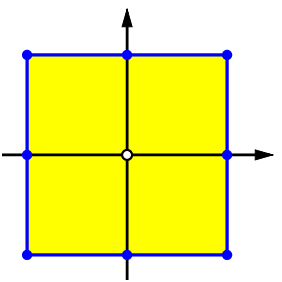} &&
%fig2dev -Leps -m0.2 Square.fig Square.eps
%fig2dev -Leps -m0.2 SquareVee.fig SquareVee.eps
      \includegraphics{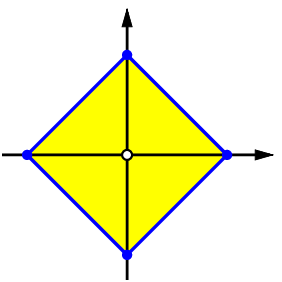} \\
      The polar of the square & \ \ldots\  & is the diamond.
   \end{tabular}
\item \begin{tabular}{ccc}
  \begin{picture}(85,65)(-10,0)
   \put(0,0){\includegraphics{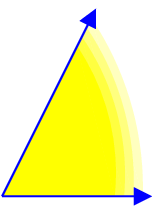}}
%fig2dev -Leps -m0.3 Cone.fig Cone.eps
   \put(-10,50){$(1,2)$}    \put(50,0){$(1,0)$} 
  \end{picture}
&&
 \begin{picture}(130,75)(-35,0)
   \put(0,0){\includegraphics{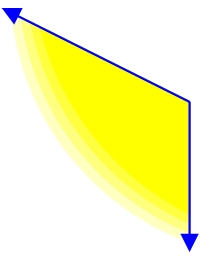}}
   \put(-35,60){$(-2,1)$}    \put(62,5){$(0,-1)$} 
%fig2dev -Leps -m0.3 ConeVee.fig ConeVee.eps
  \end{picture}
    \\
    The polar of a cone $\calK$ &\ \ldots\ & is the cone
    $\calK^\vee := \{ x \mid \langle x, y \rangle \le 0 \text{ for all }
    y \in \calK \}$ .
  \end{tabular}\smallskip

\item Suppose that $\calP$ is a polytope  whose interior contains the orign
  and $\calF$ is a face of $\calP$. Then 
  \begin{tabular}{ccc}
 \begin{picture}(100,94)
   \put(0,0){\includegraphics{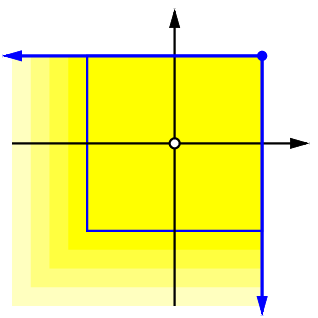}}
 %fig2dev -Leps -m0.2 TangentCone.fig TangentCone.eps
  \put(82,73){$\calF$}
 \end{picture} 
&&
 \begin{picture}(90,94)
   \put(0,10){\includegraphics{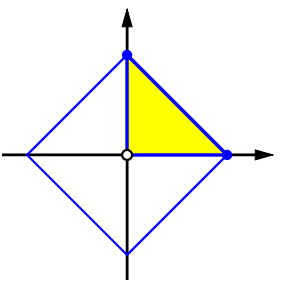}}
 %fig2dev -Leps -m0.2 TangentConeVee.fig TangentConeVee.eps
  \put(54,68){$\calF^\vee$}
 \end{picture}
\\
    the polar of the tangent cone $\calK_\calF$ & \ldots &
    is the convex hull of the origin \\
    && together with the dual face \\ &&
    $\calF^\vee := \{ x \in \calP^\vee \mid \langle x, y \rangle = 1
    \}$,\\
    && which is a pyramid over $\calF^\vee$.%
  \end{tabular}\smallskip

For this last remark, note that if $x \in \calF^\vee$
and $y \in \calK_\calF$, then    
$\langle x,y \rangle \le \langle \calF^\vee, \calF \rangle = 1$. 
Conversely, if $x \in \calK_\calF^\vee$, then 
$\langle x, \mbox{{\Large .}} \rangle$ is maximized over $\calK_\calF$ at
$\calF$ by example (2), and it is at most $1$ there.
  \end{enumerate}

In these examples, the polar of the polar is the original polyhedron. This
happens if and only if the original polyhedron contains the origin.
\begin{enumerate} \setcounter{enumi}{3}
\item The polar of the interval $[1,2]$ is the interval $[0,1/2]$, but the
  polar of $[0,1/2]$ is $[0,2]$.
\end{enumerate}
Now, we come to the main theorem of this section.

\begin{theorem}[Lawrence~\cite{La88}] \label{thm:polarity}
  The assignment $[\calP] \mapsto [\calP^\vee]$ defines a valuation.
\end{theorem}

This innocent-looking result has powerful consequences.
Suppose that $\calP$ is a polytope whose interior contains the orign. 
Then we can cover $\calP^\vee$ by pyramids $\conv(0,\calF^\vee)$ over
the codimension-one faces $\calF^\vee$ of $\calP^\vee$. The indicator
functions of $\calP$ and the cover differ by indicator functions of
pyramids of smaller dimension.
 \begin{equation} \label{eq:BGpolar}
  [\calP^\vee]\ =\ \sum_{\calF^\vee} [\conv(0,\calF^\vee)] \pm \text{
    lower dimensional pyramids} .
 \end{equation}
The Euler--Poincar\'e formula for general polytopes organizes this
inclusion-exclusion, giving the exact expression
\[
   [\calP^\vee]\ =\ \sum (-1)^{\codim \calF^\vee+1}
    [\conv(0,\calF^\vee)]\,.
\]
We illustrate this when $\calP$ is the square.
\[
  \begin{picture}(80,80)
 \put(0,7.5){\includegraphics{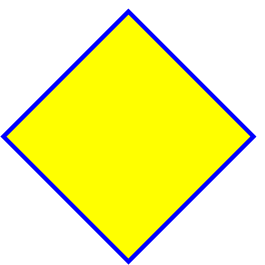}}
%fig2dev -Leps -m0.35 Diamond.fig Diamond.eps
 \end{picture}
 \quad\raisebox{42pt}{$=$}\qquad
\begin{picture}(90,90)
 \put(0,0){\includegraphics{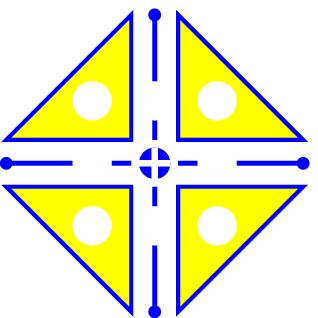}}
  %fig2dev -Leps -m0.25 DiamondDecomp.fig DiamondDecomp.eps
  \put(22,59.6){\Red{$+$}}\put(40.3,59.6){\Red{$-$}}\put(58.4,59.6){\Red{$+$}}
  \put(22,41.9){\Red{$-$}}\put(40.3,41.9){\Red{$+$}}\put(58.4,41.9){\Red{$-$}} 
  \put(22,23.8){\Red{$+$}}\put(40.3,23.8){\Red{$-$}}\put(58.4,23.8){\Red{$+$}}
 \end{picture}
%\begin{picture}(125,125)
% \put(0,0){\includegraphics{figures/DiamondDecompB}}
%  %fig2dev -Leps -m0.35 DiamondDecomp.fig DiamondDecompB.eps
%  \put(32.5,85  ){\Red{$+$}} \put(83.2,85){\Red{$+$}}
%  \put(32.5,34.2){\Red{$+$}} \put(83.2,34.2){\Red{$+$}}
%  \put(32.5,59.5){\Red{$-$}} \put(58,34.2){\Red{$-$}} 
%  \put(83.2,59.5){\Red{$-$}} \put(58,85){\Red{$-$}} 
%          \put(58,59.5){\Red{$+$}}
% \end{picture}
\]
\begin{eqnarray*}
  \includegraphics{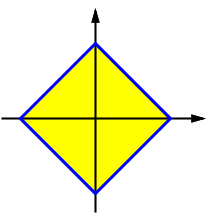}
% makeDecomp.sh
 &\raisebox{24.4pt}{=}&
  \includegraphics{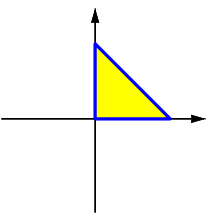}\ \raisebox{24.4pt}{$+$}\ 
  \includegraphics{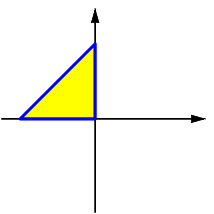}\ \raisebox{24.4pt}{$+$}\ 
  \includegraphics{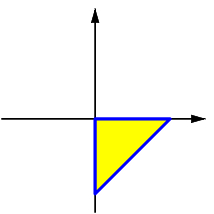}\ \raisebox{24.4pt}{$+$}\ 
  \includegraphics{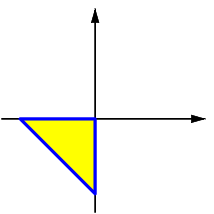}
\\
  &\raisebox{24.4pt}{$-$}&
  \includegraphics{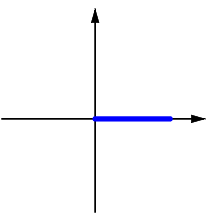}\ \raisebox{24.4pt}{$-$}\ 
  \includegraphics{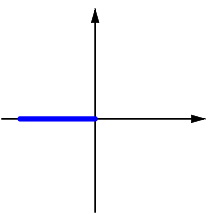}\ \raisebox{24.4pt}{$-$}\ 
  \includegraphics{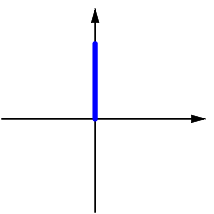}\ \raisebox{24.4pt}{$-$}\ 
  \includegraphics{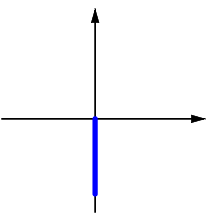}\ \raisebox{24.4pt}{$+$}\ 
  \includegraphics{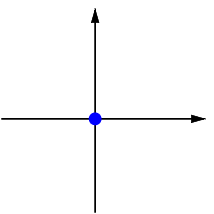}\ \raisebox{24.4pt}{$.$}
\end{eqnarray*}

If we apply polarity to \eqref{eq:BGpolar}, we get the Brianchon--Gram
Theorem~\cite{Br1837,Gr1874}.
\begin{equation} \label{eq:BG}
  [\calP]\ =\ \sum_{v \text{ vertex}} [\calK_v] \pm \text{ tangent cones
    of faces of positive dimension} .
\end{equation}
This is essentially the indicator function version of
Theorem~\ref{T:simplex}, but for general polytopes.
If we now apply the valuation $\sigma$, and recall that $\sigma$
evaluates to zero on cones that are not strictly convex, we obtain Brion's
Formula.

Next, suppose that we are given a generic direction vector $\xi$. 
On a face $\calF$ of $\calP$, the dot product with $\xi$ achieves its
maximum at a vertex $v_\xi(\calF)$.
For a vertex $v$ of $\calP$, we set
 \[
  \calF^\vee_\xi(v)\ :=\ \bigcup_{\calF \colon v_\xi(\calF) = v}
  \relint \calF^\vee .
 \]
 (The relative interior, $\relint(\calP)$, of a polyhedron $\calP$ is the
  topological interior when considered as a subspace of its affine hull.)
In words, we attach the relative interior of a low-dimensional
pyramid $\conv(0,\calF^\vee)$ to the full-dimensional pyramid
$\conv(0,v^\vee)$ which we see when we look in the $\xi$-direction from
$\conv(0,\calF^\vee)$.
In this way, we obtain an honest decomposition 
 \begin{equation} \label{eq:LVpolar}
  [\calP^\vee]\ =\ \sum_{v} \ [\conv(0,\calF^\vee_\xi(v))]\,.
 \end{equation}
For the polar of the square, this is
\begin{center}
  \includegraphics{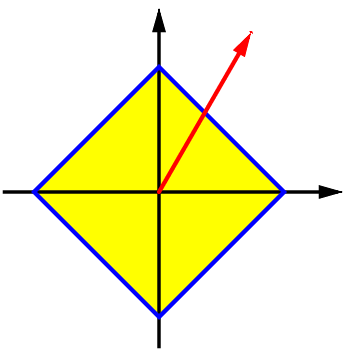}\qquad\raisebox{35pt}{$=$}\qquad 
  \includegraphics{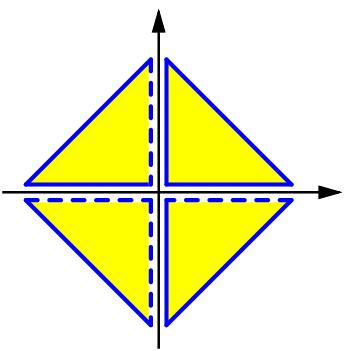}
\\
  \includegraphics{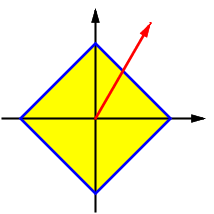}\quad \raisebox{24.4pt}{$=$}\quad 
  \includegraphics{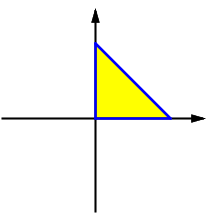}\ \raisebox{24.4pt}{$+$}\ 
  \includegraphics{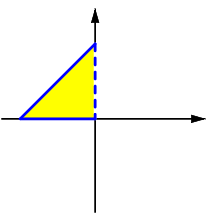}\ \raisebox{24.4pt}{$+$}\ 
  \includegraphics{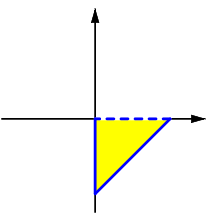}\ \raisebox{24.4pt}{$+$}\ 
  \includegraphics{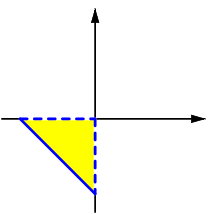}\ \raisebox{24.4pt}{$.$}\rule{0pt}{75pt}
\end{center}
To compute the polar of the half-open polyhedron
$\conv(0,\calF^\vee_\xi(v))$, we have to write its indicator function
$[\conv(0,\calF^\vee_\xi(v))]$ as a linear combination of 
indicator functions of (closed) polyhedra.
If $\calP$ is a simple polytope, then all the dual faces $\calF^\vee$
are simplices.
 It turns out that the polar of $\conv(0,\calF^\vee_\xi(v))$
is precisely the forward tangent cone $\calK_{\xi,v}$ at the vertex
$v$. So the Lawrence--Varchenko formula is just the
polar of~\eqref{eq:LVpolar}.

This gives a fairly general principle to construct Brion-type
formulas: Choose a decomposition of (the indicator function of)
$\calP^\vee$, and then polarize. 
We invite the reader to set up their own equations this way.

%%%%%%%%%%%%%%%%%%%%%%%%%%%%%%%%%%%%%%%%%%%%%%%%%%%%%%%%%%%%%%%%%%%%%
\section*{An Application}

Brion's Formula shows that 
certain data of a \emph{polytope}---the list of its integer points
encoded in a generating function---can be reduced to \emph{cones}.
We have already seen how to construct the generating function 
$\sigma_\calK (x)$ for a simple cone $\calK$. 
General cones can be composed from simple ones via triangulation and
either irrational decomposition or inclusion-exclusion.
Given a rational polytope $\calP$, Brion's Formula allows us to write 
the possibly huge polynomial $\sigma_\calP (x)$ as a sum of rational
functions, which stem from (triangulations of) the vertex cones.
{\it A priori} it is not clear that this rational-function representation
of $\sigma_\calP (x)$ is any shorter than the original polynomial.
That this is indeed possible is due to the \emph{signed decomposition}
theorem of Barvinok \cite{barvinok}.

To state Barvinok's Theorem, we call a rational $d$-cone 
$\calK = \v + \sum_{ i=1 }^d \R_{\geq 0} \w_i$
\emph{unimodular} if $\w_1, \dots, \w_d \in \Z^d$ generate the 
integer lattice $\Z^d$.
The significance of a unimodular cone $\calK$ for us is 
that its fundamental (half-open) parallelepiped contains 
precisely one integer point $\p$, and so 
the generating function of $\calK$ has a very simple and short form
\[
  \sigma_\calK (x)\ =\ \frac{ x^\p }{ \left( 1-x^{ \w_1 } \right) 
   \cdots \left( 1-x^{ \w_d } \right) } \ .
\]
In fact, the description length of this is proportional to the
description of the cone $\calK$.

\begin{theorem}[Barvinok]\label{T:Barv2} For fixed dimension $d$,
   the generating 
   function $\sigma_\calK$ for any rational cone $\calK$ in $\R^d$ can be decomposed into
   generating functions of unimodular cones in polynomial time; that is, there
   is a polynomial-time algorithm and (polynomially many) unimodular cones
   $\calK_j$ such that $\sigma_\calK(x) = \sum_j \epsilon_j \sigma_{ \calK_j
   }(x)$, where $\epsilon_j \in \{ \pm 1 \}$.  
\end{theorem}

Here \emph{polynomial time} refers to the input data of $\calK$, 
that is, the algorithm runs in time polynomial in the input length 
of, say, the halfspace description of $\calK$.

Brion's Formula implies  that an identical complexity 
statement can be made about the generating function $\sigma_\calP(x)$
for any rational polytope $\calP$.
{}From here it is a short step (which nevertheless needs some
justification) to see that one can \emph{count} integer points 
in a rational polytope in polynomial time.

We illustrate Barvinok's short signed decomposition for the cone
$\calK:=(0,0)+\R_{\ge0}(1,0) + \R_{\ge0}(1,4)$, ignoring cones of
smaller dimension.
\[
 \begin{picture}(95,128)(-6,-12)
  \put(0,0){\includegraphics{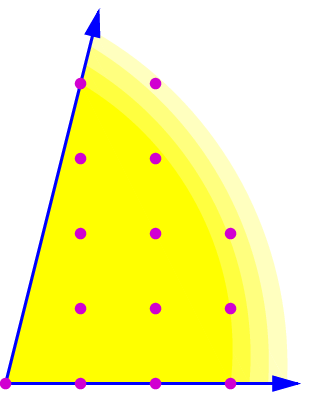}}
  \put(-7,87){$(1,4)$} \put(55,-10){$(3,0)$}
 \end{picture}
 \quad\raisebox{55 pt}{$=$}\quad 
  \begin{picture}(130,128)(-26,-12)
   \put(0,0){\includegraphics{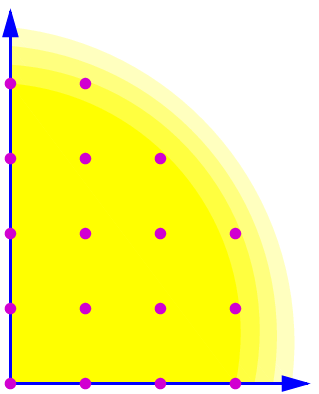}}
  \put(-28,87){$(0,4)$} \put(55,-10){$(3,0)$}
 \end{picture}
  \ \raisebox{55 pt}{$-$}\ 
  \begin{picture}(70,128)(-26,-12)
   \put(0,1.5){\includegraphics{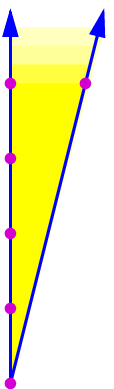}}
  \put(-28,87){$(0,4)$}  \put(27,87){$(1,4)$} 
 \end{picture}
 \]
While $\calK$ is the difference of two unimodular cones,
it has a unique decomposition as a sum of {\it four} unimodular cones.  
\[
 \begin{picture}(95,128)(-6,-12)
  \put(0,0){\includegraphics{figures/Cone4}}
  \put(-7,87){$(1,4)$} \put(55,-10){$(3,0)$}
 \end{picture}
 \quad\raisebox{55 pt}{$=$}\quad 
  \begin{picture}(88,128)(0,-12)
   \put(0,0){\includegraphics{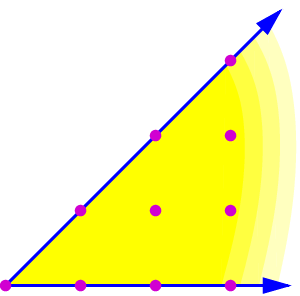}}
     \put(35,70){$(3,3)$} \put(55,-10){$(3,0)$}
 \end{picture}
 \ \raisebox{55 pt}{$+$}\ 
  \begin{picture}(80,128)(5,-12)
   \put(0,1.5){\includegraphics{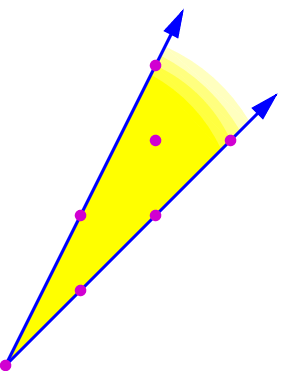}}
     \put(15,87){$(2,4)$} \put(60,53){$(3,3)$}
 \end{picture}
 \ \raisebox{55 pt}{$+$}\ 
  \begin{picture}(58,128)(0,-12)
   \put(0,1.5){\includegraphics{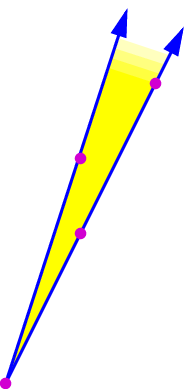}}
    \put(-7,65){$(1,3)$}   \put(27,42){$(1,2)$} 
 \end{picture}
 \ \raisebox{55 pt}{$+$}\ 
   \begin{picture}(47,128)(3,-12)
   \put(0,1.5){\includegraphics{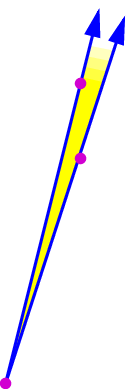}} 
   \put(-7,87){$(1,4)$}    \put(27,60){$(1,3)$}  
 \end{picture}
\]
In general the cone $(0,0)+\R_{\ge0}(1,0) + \R_{\ge0}(1,n)$ is the
difference of two unimodular cones, but it has a unique decomposition
into $n$ unimodular cones.

Arguably the most famous consequence of Barvinok's Theorem applies
to \emph{Ehrhart quasi\-polynomials}---the counting functions
$L_\calP (t) := \# \left( t \calP \cap \Z^d \right)$ in the 
positve-integer variable $t$ for a given rational polytope \cite{beckrobins} 
$\calP$. One can show that the generating 
function $\sum_{ t \ge 1 } L_\calP (t) \, x^t$ is a rational 
function, and Barvinok's Theorem implies that this rational 
function can be computed in polynomial time.
Barvinok's algorithm has been implemented in the software packages 
{\tt barvinok} \cite{verdoolaege} and {\tt LattE} \cite{lattemanual}.
The method of irrational decomposition has also been implemented in 
{\tt LattE}, considerably improving its performance~\cite{Koeppe}.

%%%%%%%%%%%%%%%%%%%%%%%%%%%%%%%%%%%%%%%%%%%%%%%%%%%%%%%%%%%%%%%%%%%%%%%%%%%%%%
%
\section*{Acknowledgments}
%
%%%%%%%%%%%%%%%%%%%%%%%%%%%%%%%%%%%%%%%%%%%%%%%%%%%%%%%%%%%%%%%%%%%%%%%%%%%%%%
Research of Haase supported in part by NSF grant DMS-0200740
and DFG Emmy Noether fellowship.
Research of Sottile supported in part by the Clay Mathematical
nstitute and NSF CAREER grant DMS-0538734.

%%%%%%%%%%%%%%%%%%%%%%%%%%%%%%%%%%%%%%%%%%%%%%%%%%%%%%%%%%%%%%%%%%%%%%%%%%%%%%

\def\cprime{$'$}
\providecommand{\bysame}{\leavevmode\hbox to3em{\hrulefill}\thinspace}
\providecommand{\MR}{\relax\ifhmode\unskip\space\fi MR }
% \MRhref is called by the amsart/book/proc definition of \MR.
\providecommand{\MRhref}[2]{%
  \href{http://www.ams.org/mathscinet-getitem?mr=#1}{#2}
}
\providecommand{\href}[2]{#2}

\end{document}